\documentclass{amsart}
\usepackage{amssymb} 
\usepackage{amscd}

\newtheorem{thm}{Theorem}

\newtheorem{lemma}{Lemma}
\newtheorem{remark}{\it Remark}

\newcommand{\R}{\mathbb{ R}}

\newcommand{\D}{\mathfrak{ d}}

\DeclareMathOperator{\pr}{pr}
\DeclareMathOperator{\Span}{span}
\DeclareMathOperator{\diag}{diag}

\title[Cartan Model of the Canonical Vector Bundles]{On the Cartan Model of the Canonical Vector Bundles over
Grassmannians}

\author{Bo\v zidar Jovanovi\' c}

\address{Mathematical Institute SANU, Kneza Mihaila 35, 11000 Belgrade, Serbia and Montenegro}
\email{bozaj@mi.sanu.ac.yu}

\begin{document}

\maketitle

\begin{abstract}
We give a representation of canonical vector bundles $\mathcal
C_{n,p}$ over Grassmannian manifolds $G(n,p)$ as non-compact
affine symmetric spaces as well as their Cartan model in the group
of the Euclidean motions $SE(n)$.
\end{abstract}

\section{Introduction}

The Cartan model of Grassmannian manifolds $G(n,p)$ in the special
orthogonal group $SO(n)$ is well known. Remarkably, we find that
there is a representation of the canonical vector bundles
$\mathcal C_{n,p}$ over $G(n,p)$ as symmetric spaces, namely
$\mathcal C_{n,p}=SE(n)/S(O(p)\times O(n-p)) \otimes_s \R^{n-p}$
and Cartan model realization in the group of  Euclidean motions
$SE(n)$. To the author knowledge, this interesting fact is not
observed yet (e.g., see \cite{KN, H}).

The different homogeneous space
representation of the canonical line bundles over projective spaces
can be found in \cite{G}. The Cartan-type model of the M\"oebius
strip in $SE(2)$ is recently obtained in \cite{FZ}. Note that, due
to \cite{YK}, tangent bundles of Grassmannians have natural affine
symmetric space structures.

\section{Grassmannian Varieties}

The points of the Grassmannian variety $G(n,p)$ are by
definition $p$-dimensional planes $\pi$ passing
through the origin of $\R^n$.
In particular, for $p=1$, we have the projective space
$\mathbb{RP}^{n-1}$, the set of lines
through the origin in $\R^n$.

Grassmannian manifolds are basic examples of compact symmetric
spaces. The usual action of the group $SO(n)$ on $\R^n$ yields a
transitive action on the set of all $p$-dimensional planes, i.e.,
on $G(n,p)$. Let
$$
E_1=(1,0,\dots,0)^T, \quad E_n=(0,\dots,0,1)^T.
$$

Take the plane $\pi_0=\Span\{E_1,\dots,E_p\}$. Then the
isotropy group of $\pi_0$ consists of matrixes
$$
\begin{pmatrix}
A & 0 \\
0 & B
\end{pmatrix}, \quad A\in O(p), \quad B\in O(q), \quad  \det A\cdot\det B=1.
$$
It follows that $G(n,p)\cong SO(n)/(S(O(p)\times O(q))$.
Further, let
$$
J_{p,q}=\begin{pmatrix}
-\mathbb{I}_p & 0 \\
0 & \mathbb{I}_q
\end{pmatrix},
$$
where $\mathbb{I}_k=\diag(1,1,\dots,1)$. Then $\sigma_0: SO(n)\to SO(n)$,
$$
\sigma_0(R)= J_{p,q} R J_{p,q}
$$
is an involutive automorphism with $SO(n)^{\sigma_0}=S(O(p)\times O(q))$
and the triple $(SO(n),S(O(p)\times O(q)),\sigma_0)$ is a symmetric space.

\subsection*{Cartan model of Grassmannians}
Let
$$
\D_p^0=\Span\{E_i\wedge E_j\, \vert\, 1\le i\le p <j \le n=p+q\}\subset so(n).
$$
Then $so(n)=so(p)+so(q)+\D_p^0$ is the symmetric pair decomposition of
the Lie algebra $so(n)$ on $(+1)$ and $(-1)$ eigenspaces of $d\sigma_0$
at the identity $\mathbb{I}_n$.

Consider the $\sigma_0$-twisted conjugation action
$A\bullet R=A R \sigma_0 (A)^{-1}$, $R,A\in SO(n)$.
Let
$$
\mathcal Q_0=\{R\in SO(n)\, \vert\, \sigma_0(R)=R^{-1}\}=
\{R\in SO(n)\, \vert\, (RJ_{p,q})^2=\mathbb{I}_n\}.
$$
It can be easily verified that $\mathcal Q_0$ is invariant under
the $\sigma_0$-twisted action.

The orbit through identity
$$
\mathcal S_p^0=SO(n)\bullet \mathbb{I}_n=
\{A\sigma(A)^{-1}=A J_{p,q} A^{-1} J_{p,q} \,\vert\, A\in SO(n)\},
$$
is isomorphic to $G(n,p)$ as a $SO(n)$-space, relative to the
$\sigma_0$-twisted conjugation action and $\mathcal S_p^0$
coincides with the identity connected component of $\mathcal Q_0$
({\it the Cartan model of a symmetric space}, e.g., see
\cite{Fo}). Furthermore, $\mathcal S_p^0$ is equal to the image of
$\D_p^0$ under the exponential mapping.

Take the translation $\mathcal S_p^0 \cdot J_{p,q}=\{A J_{p,q}
A^{-1} \,\vert\, A\in SO(n)\}$. The matrix $AJ_{p,q}A^{-1}$ is
symmetric and has $(-1)$ eigenvalue on the plane $\pi=A\cdot
\pi_0=\Span\{A\cdot E_1,\dots,A\cdot E_p\}$. Thus, the
diffeomorphism $\rho_0: \mathcal S_p^0 \to G(n,p)$ can be seen as
follows:
$$
\rho_0(R)=\pi,
$$
where $\pi$ is the unique plane satisfying
$RJ_{p,q}(X)=-X$, $X\in \pi.$

\subsection*{Projective Spaces}
For $p=1$, $J_{1,q}$ is the reflection $S_{1}$ with respect to the
plane orthogonal to $E_1$. Further, the elements of $\D_1^0$ can
be taken to be of the form $-\theta E_1\wedge U$, $\vert
U\vert=1$, $U\perp E_1$. Then $R_{\theta,U}=\exp(-\theta E_1\wedge
U)$ is the rotation in the plane spanned by $E_1$ and $U$:
\begin{equation}
R_{\theta,U}(E_1)=\cos\theta E_1+\sin\theta U, \quad R_{\theta,U}(U)=-\sin\theta E_1+\cos\theta U,
\label{rotation}
\end{equation}
which fix the orthogonal complement to $\Span\{E_1,U\}$. The
rotation can be represented as a composition: $R_{\theta,U}=S_2
\circ S_1$, where $S_2$ is the reflection with respect to the
plane orthogonal to the vector $V=\cos\frac\theta2 E_1+
\sin\frac\theta2 U$. Since $R_{\theta,U}J_{1,q}=S_2\circ S_1\circ
S_1=S_2$ and $S_2(V)=-V$, we get
\begin{equation}
\rho_0(R_{\theta,U})= [V]=\left[\cos\frac\theta2 E_1+
\sin\frac\theta2 U\right]. \label{Cartan_projective}
\end{equation}
Here $[V]$ denotes the line $\{\mu V \, \vert \,\mu\in \R\}$.

\section{Cartan Model of the Canonical Vector Bundles}

Consider $SE(n)$, the Lie group of the motions in the Euclidean space
$(\mathbb{R}^n,\langle\cdot,\cdot\rangle)$.
It is a semi-direct product
of the special orthogonal group $SO(n)$ (rotations) and the abelian group $\mathbb{R}^n$
(translations) $SE(n)=SO(n)\otimes_s \R^n$. We use the following usual matrix notation
for the elements $g\in SE(n)$:
$$
g=(R,X)=\begin{pmatrix}  R & X \\
                     {0} & 1 \\
               \end{pmatrix}, \quad R\in SO(n), \quad X\in\R^n.
$$
The Lie algebra $se(n)=so(n)\oplus_s \R^n$ consist of the $(n+1)\times(n+1)$ matrixes
$$
\xi=(\omega,v)=\begin{pmatrix}   \omega & v \\
                     {0} & 0  \\
               \end{pmatrix}, \quad \omega\in so(n), \quad v\in\R^n.
$$
The group multiplication and Lie bracket correspond to the usual multiplication and
Lie bracket for the matrixes:
\begin{eqnarray*}
&&(R_1,X_1)\cdot (R_2,X_2)=(R_1R_2, X_1+R_1 X_2), \\
&&[(\omega_1,v_1),(\omega_2,v_2)]=([\omega_1,\omega_2],\omega_1 v_2-\omega_2 v_1).
\end{eqnarray*}

\begin{lemma}
The mapping $\sigma: \, SE(n)\to SE(n)$ given by
\begin{equation}
\sigma((R,X))=(\sigma_0(r),J_{p,q} X)=(J_{p,q} R J_{p,q}, J_{p,q} X).
\label{involution}
\end{equation}
is an involutive automorphism and the
set of fixed point consist of matrixes of the form
$$
\begin{pmatrix}
A & 0 & 0\\
0 & B & X\\
0 & 0 & 1
\end{pmatrix}, \quad A\in O(p), \quad B\in O(q), \quad  \det A\cdot\det B=1, \quad X\in \R^q,
$$
i.e., $SE(n)^\sigma=S((O(p)\times O(q))\otimes_s \R^q$
\end{lemma}

Therefore, the triple $(SE(n),S((O(p)\times O(q))\otimes_s
\R^q,\sigma)$ is a non-compact affine symmetric space (we follow
the notation of \cite{KN}). The differential $d\sigma$ at the
identity $(\mathbb{I}_n,0)$ is an involutive
automorphism of the Lie algebra $se(n)$. We have symmetric pair
decomposition of $se(n)$ on its $(+1)$ eigenspace (the Lie algebra
of $SE(n)^\sigma$) and $(-1)$ eigenspace:
$$
\D_p=\Span\{(E_i\wedge E_j,E_k)\, \vert\, 1\le i\le p <j \le n=p+q, 1\le k \le p\}
\cong \D_p^0\oplus \R^p.
$$

Let
$$\mathcal Q=\{g=(R,Y)\in SE(n)\, \vert\, \sigma(g)=g^{-1}\}.$$

The set $\mathcal Q$ is preserved under the
$\sigma$-twisted conjugation action:
\begin{eqnarray*}
(A,X)\bullet (R,Y) &=& (A,X)\cdot (R,Y)\cdot \sigma((A,X)^{-1})\\
&=& (ARJ_{p,q} A^{-1} J_{p,q}, X+AY-ARJ_{p,q}A^{-1}X).
\end{eqnarray*}

The Cartan model of symmetric spaces  is usually given for
reductive Lie groups. Similarly we have

\begin{thm} {\rm(The Cartan Model)}
The orbit through the identity
$$
\mathcal S_p=SE(n)\bullet (\mathbb{I}_n,0)=
\{(A J_{p,q} A^{-1} J_{p,q},X-A J_{p,q} A^{-1} X) \,\vert\, (A,X)\in SE(n)\}.
$$
is isomorphic to $SE(n)/SE(n)^\sigma$ as a $SE(n)$-space, relative
to the $\sigma$-twisted conjugation action. Furthermore, $\mathcal
S_p$ is equal to the identity component of $\mathcal Q$ and it is
equal to the image of $\D_p$ under the exponential mapping.
\end{thm}

\begin{lemma} The exponential mapping $\exp: se(n)\to SE(n)$ is surjective.
\end{lemma}
\noindent{\it Proof.}
A simple computation shows
$$
\xi^m=(\omega,v)^m=(\omega^m,\omega^{m-1} v), \quad m\in\mathbb{N}.
$$
Therefore
$$
\exp(\xi)=(\exp(\omega),Y),
$$
where the vector $Y=Y_\omega(v)$ is equal to
\begin{equation}
Y_\omega(v)=v+\frac12 \omega v+\frac{1}{3!} \omega^2 v +\dots+ \frac{1}{m!} \omega^{m-1} v+ \dots.
\label{0}
\end{equation}

Since $\exp: so(n)\to SO(n)$ is surjective, we only need to prove
that the linear mapping (\ref{0}), for the fixed $R\in SO(n)$ and
properly chosen $\omega$, $R=\exp(\omega)$, has the maximal rank.

Let $e_1,\dots,e_n$ be the orthonormal base , in which the matrix $R$ has the
canonical form
$$
R=\diag(R(\theta_1),R(\theta_2),\dots,R(\theta_k),1,1,\dots,1),
$$
where $R(\theta_i)$ are rotations in the planes $\Span\{e_{2i-1},e_{2i}\}$:
$$
R(\theta_i)=\begin{pmatrix}
\cos \theta_i& -\sin\theta_i \\
\sin\theta_i & \cos\theta_i
\end{pmatrix}, \quad \vert\theta_i\vert < 2\pi, \quad i=1,\dots,k.
$$
Then we can take
$\omega=\diag(\Pi(\theta_1),\Pi(\theta_2),\dots,\Pi(\theta_k),0,0,\dots,0),$
where
$$
\Pi(\theta_i)=\begin{pmatrix}
0 & -\theta_i \\
\theta_i & 0
\end{pmatrix}, \quad i=1,\dots,k.
$$

Let $Y_i(v)=\langle Y_\omega(v),e_i\rangle$.
For a given $v=v_1e_1+\dots v_n e_n$ we have
\begin{equation}
Y_i(v)=v_i, \quad i=2k+1,\dots,n.
\label{2}
\end{equation}
Further, from (\ref{0}) we get that $Y_\omega(v)$ satisfies the relation
\begin{equation}
\omega Y_\omega(v)=(\exp(\omega)-\mathbb{I}_n) v,
\label{1}
\end{equation}
or, in coordinates:
\begin{eqnarray*}
&&-\theta_i Y_{2i}(v)=\cos\theta_i v_{2i-1}-\sin\theta_i v_{2i}-v_{2i-1} \\
&&\theta_i Y_{2i-1}(v)=\sin\theta_i v_{2i-1}+\cos\theta_i v_{2i}-v_{2i}, \quad i=1,\dots,k.
\end{eqnarray*}
By using the trigonometric identities
$\sin\theta=2\sin\frac{\theta}2\cos\frac{\theta}2$,
$1-\cos\theta=2\sin^2\frac{\theta}2$,
we can write the components of the vector $Y$ in the compact form:
\begin{eqnarray}
&&Y_{2i-1}(v)=2\frac{\sin\frac{\theta_i}2}{\theta_i}
\left(\cos\frac{\theta_i}2 v_{2i-1}-\sin\frac{\theta_i}2 v_{2i}\right) \label{3}\\
&&Y_{2i}(v)=2\frac{\sin\frac{\theta_i}2}{\theta_i}\left(\sin\frac{\theta_i}2 v_{2i-1}+\cos\frac{\theta_i}2 v_{2i}\right), 
\quad i=1,\dots,k.\label{4}
\end{eqnarray}

From (\ref{2}), (\ref{3}), (\ref{4}), it follows
that $Y_\omega$ has no kernel. $\Box$

\medskip

\noindent{\it Proof of the Theorem.} In proving the Theorem, we
mainly follow standard arguments given for compact (or reductive)
Lie groups (e.g., see \cite{Fo}).

{\bf (i)} Let $\tau: SE(n)\to \mathcal S_p$ be the mapping defined
by $\tau(g)=g\sigma(g^{-1})$. It is clear that $\tau$ is constant
on left cosets modulo $SE(n)^\sigma$ ($\tau(g_1)=\tau(g_2)$ if and
only if $\sigma(g_1 g_2^{-1})=g_1 g_2^{-1}$, i.e, $g_1 g_2^{-1}\in
SE(n)^\sigma$) and that the induced morphism $\hat\tau:
SE(n)/SE(n)^\sigma\to\mathcal S_p$ is bijective and satisfies
$$
\hat \tau(g_1\cdot g_2 SE(n)^\sigma)=g_1 \bullet  \hat\tau (g_2 SE(n)^\sigma).
$$
Further, $\hat \tau$ is a diffeomorphism
from the dimensional reasons. (It can be easily seen that the tangent space of $\mathcal S_p$
at the identity of the group is $\D_p$ so the differential $d\hat\tau\vert_{SE(n)^\sigma}$
is surjective.)

{\bf (ii)} Suppose that $(R,Y)$ belongs to the identity component of $\mathcal Q$. Then
$$
\sigma(R,Y)=(R^{-1},-R^{-1}Y), \quad {\rm{i.e.,}}\quad
\sigma_0(R)=R^{-1} \quad {\rm{and}} \quad J_{p,q}Y=-R^{-1}Y.
$$

Recall that $\mathcal S_p^0$ coincides with the identity component
of $\mathcal Q_0$. Therefore $R=A J_{p,q} A^{-1} J_{p,q}$, for
some $A\in SO(n)$. Then the condition $J_{p,q}Y=-R^{-1}Y$ is the
same as
$$
J_{p,q}(Y+J_{p,q}R^{-1}Y)=J_{p,q}(Y+A J_{p,q} A^{-1} Y)=0,
$$
which leads $Y\in \pi=\rho_0(R)=A\cdot \pi_0$.
On the other side, for $X\in\R^n$,
we have
\begin{equation}
X-A J_{p,q} A^{-1} X=2\pr_\pi X,
\label{projection}
\end{equation}
where $\pr_\pi$ denotes the orthogonal projection
to $\pi$.
Therefore, $(R,Y)\in\mathcal S_p$\footnote{The alternative proof of this statement is to show that
the tangent space to a $\sigma$-twisted $SE(n)$-orbit coincides the tangent space to $Q$}.
The another inclusion is trivial:
$g=g'\sigma(g'^{-1})$ implies $\sigma (g)=\sigma (g') g'^{-1}=g^{-1}$.

{\bf (iii)} First, we shall prove the inclusion
$\exp(\D_p)\subset\mathcal S_p$. Let $g=\exp(\xi)$, $\xi\in\D_p$.
Consider the element $g'=\exp(\xi/2)$. Then
$$
\tau(g')=\exp(\xi/2)
\sigma(\exp(-\xi/2))=\exp(\xi/2)\exp(\xi/2)=(g')^2=g,
$$
that is $g\in \mathcal S_p$. Here we used the identity
$\sigma(\exp(\xi))=\exp(d\sigma\vert_{(\mathbb{I}_n,0)}\xi)$.

Now, let $R$ be an arbitrary element in $\mathcal S^0_p$. From
$\mathcal S_p^0=\exp(\D_p^0)$ and Lemma 2, for a properly chosen
$\omega\in \D_p^0$, $R=\exp(\omega)$, we have that the linear
mapping (\ref{0}) define an isomorphism between
$\Span\{E_1,\dots,E_p\}$ and $\pi=\rho_0(R)$. Therefore $\exp: \D_p\to
\mathcal S_p$ is a surjective map. $\Box$

\medskip

Recall that the {\it canonical vector bundle} $\mathcal C_{n,p}$
over $G(n,p)$ at the point $\pi\in G(n,p)$
has the fibre equal to $\pi$, now considered as a vector space:
$$
\mathcal C_{n,p}=\{(\pi,X)\in G(n,p)\times\R^n\, \vert\, X\in\pi\}.
$$

\begin{lemma}
The variety $\mathcal S_p$ is diffeomorphic to the canonical vector bundle
$\mathcal C_{n,p}$.
\end{lemma}

\noindent{\it Proof.}
According to (\ref{projection}), the mapping
$\rho: \mathcal S_p \to \mathcal C_{n,p}$, defined by
$$
\rho(R,Y)=(\rho_0(R),Y)
$$
establish the diffeomorphism between $\mathcal S_p$ and $\mathcal C_{n,p}$. $\Box$

\medskip

From the above considerations, we see that the canonical vector
bundles over Grassmannians in a canonical way can be considered as
symmetric spaces.

\begin{thm}
\begin{equation}
(A,X)\ast (\pi,Y)= (A\pi, A Y+ 2\pr_{A\pi} X), \quad Y\in\pi\subset\R^n
\label{action}\end{equation}
defines a transitive  $SE(n)$-action
on the canonical vector bundle $\mathcal C_{n,p}$ over $G(n,p)$
such that $\rho$ becomes a $SE(n)$-invariant diffeomorphism:
$$
\rho((A,X)\bullet (R,Y))=(A,X)\ast\rho(R,Y), \quad (A,X)\in SE(n), \, (R,Y)\in\mathcal S_p.
$$
Therefore, the $SE(n)$-action (\ref{action}) realizes $\mathcal
C_{n,p}$ as a non-compact affine symmetric space
$(SE(n),S((O(p)\times O(q))\otimes_s \R^q,\sigma)$ .
\end{thm}

\begin{remark}{\rm
The different homogeneous space
representation of the canonical line bundles over projective spaces
can be found in \cite{G}. Namely, $\mathbb{RP}^{n} \setminus x_0$
is diffeomorphic to $\mathcal C_{n-1,1}$, where
$x_0\in\mathbb{RP}^n$ is an arbitrary point. Then projective
transformations of $\mathbb{RP}^n$ which leave $x_0$ invariant
acts transitively on $\mathbb{RP}^n\setminus x_0\approx \mathcal
C_{n-1,1}$.}
\end{remark}

\begin{remark}
{\rm
The description of $\exp(\mathfrak d_p)$ is
important in the study of discrete nonholonomic LL systems on $SE(n)$ (see \cite{FZ}).
}\end{remark}

\subsection*{Canonical Line Bundles}
For $p=1$, there is a direct construction of a diffeomorphism
between $\exp(\D_1)$ and the canonical line bundle $\mathcal
C_{n,1}$\footnote{This example is motivated by \cite{FZ} and was
the staring point in writing this note.}. The elements of $\D_1$
can be taken to be of the form $\xi=(-\theta E_1\wedge U,\lambda
E_1)$, $\theta,\lambda\in\R$, $\vert U\vert=1$, $U\perp E_1$. Let
$(R_{\theta,U},Y)=\exp(-\theta E_1\wedge U,\lambda E_1)$. Then the
relation (\ref{1}) reads
\begin{equation}
-\theta E_1\wedge U(Y)=\lambda R_{\theta,U} E_1-\lambda E_1.
\label{Y}
\end{equation}
Therefore, taking into account (\ref{rotation}) and (\ref{Y}) we obtain
\begin{equation}
-\theta\langle U,Y \rangle E_1+\theta\langle E_1,Y\rangle U=
\lambda(\cos\theta-1)E_1+\lambda\sin\theta U.
\label{***}
\end{equation}
From (\ref{0}) we get that $Y$ belongs to $\Span\{E_1,U\}$ and
(\ref{***}) gives
\begin{equation}
Y=\lambda \frac{\sin\theta}{\theta} E_1+\lambda\frac{1-\cos\theta}{\theta} U=
\lambda \frac{2\sin\frac{\theta}{2}}{\theta}(\cos\frac{\theta}{2} E_1+\sin\frac\theta2 U).
\label{****}
\end{equation}
Finally, in the view of (\ref{Cartan_projective}) and (\ref{****})
we get $\exp(\D_1) \approx \mathcal C_{n,1}$.

\subsection*{Acknowledgments}
I am very grateful to Yuri N. Fedorov for stimulating discussions.
The research was supported by the Serbian Ministry of Science,
Project "Geometry and Topology of Manifolds and Integrable
Dynamical Systems".

\end{document}